\documentclass[11pt]{amsart}  
\usepackage{amssymb, latexsym, amsmath, amsfonts,cases}

\makeatletter
\renewcommand*\env@matrix[1][*\c@MaxMatrixCols c]{%
  \hskip -\arraycolsep
  \let\@ifnextchar\new@ifnextchar
  \array{#1}}
\makeatother

\newtheorem{thm}{Theorem}[section]

\newtheorem{prop}[thm]{Proposition}
\theoremstyle{definition}    
\newtheorem{defn}[thm]{Definition}
\newtheorem{rem}[thm]{Remark} 
 
\numberwithin{equation}{section}

\newcommand{\abs}[1]{\left\vert#1\right\vert}

\newcommand{\mbb}{\mathbb}

\begin{document}
\title{Analyzing the Wu metric on a class of
eggs in $\mathbb{C}^n$ -- I}
\thanks{Both the authors were supported by the DST-INSPIRE Fellowship of the Government of India.}
\author{G. P. Balakumar And Prachi Mahajan}

\address{G. P. Balakumar:  Indian Statistical Institute, Chennai -- 600113, India.}
\email{gpbalakumar@isichennai.res.in}

\address{Prachi Mahajan: Department of Mathematics, Indian Institute of Technology -- Bombay, Mumbai 400076, India.}
\email{prachi.mjn@iitb.ac.in}
\pagestyle{plain}

\begin{abstract}
We study the Wu metric on convex egg domains of the form 
\[
E_{2m} = \big\{ z \in \mbb C^n : \abs{z_1}^{2m} + \abs{z_2}^2 + \ldots + \abs{z_{n-1}}^2 + \abs{z_n}^{2} <1  \big\}
\]
where $m \geq 1/2, m \neq 1$. The Wu metric is shown to be real analytic everywhere except on a lower dimensional subvariety where it fails to be $C^2$-smooth. Overall however, the Wu metric is shown to be continuous when $m=1/2$ and even $C^1$-smooth for each $m>1/2$, and in all cases, a non-K\"ahler Hermitian metric with its holomorphic curvature strongly negative in the sense of currents. This gives a natural answer to a conjecture of S. Kobayashi and H. Wu for such $E_{2m}$.
\end{abstract}
\maketitle

\section{Introduction}
\noindent Connections between various forms of the Schwarz lemma and hyperbolicity in complex analysis, 
i.e., negativity of curvature of a biholomorphically invariant metric, are too well known to need any
elaboration (a brief overview can be found in \cite{R}). Nevertheless, certain fundamental conjectures towards building 
a {\it theory} of hyperbolicity pertinent to complex analysis in standard differential geometric terms,
still remain unsettled. The present article is 
one of a two-part note, concerning the following conjecture of Kobayashi \cite{Ko}
and its modified form by Wu \cite{Wu}.
\medskip\\
\noindent \underline{Conjecture (K-W)}:
On every Kobayashi complete hyperbolic complex manifold, there exists a $C^k$-smooth (for some $k \geq 0$) complete 
Hermitian metric with its holomorphic curvature \footnote{The term `holomorphic
curvature' stands precisely for the holomorphic {\it sectional} curvature and is said 
to be strongly negative, if it is bounded above by a negative constant.} bounded above by a negative constant in the sense of currents. 
 \medskip\\
\noindent Needless to say, for the above conjecture to contribute to complex analysis, the metric must transform in a `nice' manner
under holomorphic mappings. The Kobayashi metric is the largest pseudo-distance with the nice feature of being distance 
decreasing under all holomorphic mappings. However, the Kobayashi metric fails to be Hermitian in general, obstructing applicability of
standard (Hermitian) differential geometric techniques.  A  method proposed by H. Wu in \cite{Wu} 
to redress this problem is as follows. For an infinitesimal metric $\eta$, form a new metric whose unit ball about the origin, in each tangent space, is 
the ellipsoid of minimal volume containing the unit ball with respect to $\eta$. By an ellipsoid, we mean the image of a ball in the 
standard $l^2$-norm, under a linear automorphism of $\mathbb{C}^n$. As ellipsoids furnish the simplest compact convex sets with (quadratic) polynomial
boundary, Wu's methodology is perhaps the first best attempt at constructing Hermitian (invariant) metrics. In the present article, we work with $ \eta $ as the infinitesimal Kobayashi metric \footnote{We shall reserve the term `metric' for what was termed as a sub-metric in \cite{Wu} and
use the word `distance' for what is termed as a `metric' in general topology and metric geometry.} and refer the resultant metric as 
the Wu metric.
\medskip\\
\noindent  An explicit description of the Wu metric on the Thullen domains \footnote{These are also referred to as `complex ellipsoids' 
in the complex analysis literature. However, since the ellipsoid in Wu's construction is defined by a polynomial of degree two, we  
adopt a different terminology.} in dimension $n=2$, namely $ \{ (z_1, z_2) \in \mathbb{C}^2: 
|z_1|^{2m} + |z_2|^{2} < 1 \}$, where $ m > 0 $, was obtained by Cheung and Kim in \cite{CK1} and \cite{CK2}. 
This is a natural example arising in the pursuit of well-adapted metrics in complex analysis, the background for which can be found in \cite{Wu} and \cite{Kim}.
 Let us make a few 
remarks. Although the statement of the conjecture (KW) as it stands, had been affirmed for smoothly bounded Thullen domains way back in the 1980's 
owing to the fact that the curvatures of the Bergman and the K\"ahler-Einstein metric do lie trapped between a pair of negative constants (\cite{AzSu}, \cite{Bl}),
these metrics are not functorial: the Hermitian metric sought for, in the conjecture (K-W) is desired to transform holomorphic 
mappings into Lipschitz mappings with respect to the distance induced by the metric, {\it wherever such a metric is defined}. Being 
functorial and the simplest modified (Hermitian) form of the Kobayashi metric, the Wu metric is the first natural candidate to be 
investigated for (K-W). However, the Wu metric may fail to be upper semicontinuous \cite{Ju}, notwithstanding 
the fact that the Kobayashi metric is always upper semicontinuous. But then again, the Wu metric seems to reflect the boundary
geometry {\it better} in special instances. For instance, the holomorphic curvature for the Kobayashi metric is well defined 
\cite{JP} and for any bounded convex domain in $\mathbb{C}^n $, the holomorphic curvature of the Kobayashi metric is identically a negative 
constant. On the other hand, for $C^2$-smooth convex Thullen domains in $ \mathbb{C}^2 $, the Wu metric is $C^1$-smooth Hermitian whose curvature 
is a negative constant only on an appropriate neighbourhood of the strongly pseudoconvex portion of the boundary (see \cite{CK1}). The complement of 
the closure of this neighbourhood forms a approach region for a thin subset of the boundary 
consisting of weakly pseudoconvex points and the curvature of Wu metric is nowhere constant on this region.  Nevertheless, in all 
cases, the Wu metric retains the strong negativity of curvature. \medskip\\
\noindent The purpose of this note and its counterpart \cite{BM2}, is to carry forward the case study of the 
Wu metric in conjunction with (K-W) as in \cite{CK1} and \cite{CK2}, to egg domains of the form
\begin{equation} \label{Egg}
E_{2m} = \big\{ z \in \mbb C^n : \abs{z_1}^{2m} + \abs{z_2}^2 + \ldots + \abs{z_{n-1}}^2 + \abs{z_n}^{2} <1  \big\}.  
\end{equation} 
More specifically, we reinstate the viewpoint of \cite{CK1}, \cite{CK2}: being functorial and {\it respecting Kobayashi hyperbolicity}, the Wu metric is conceivably the most natural answer to (K-W) 
when restricted to a suitable class of domains -- the class for this article being the eggs $E_{2m}$ as above, with $m \geq 1/2$. 
This condition on $m$ is precisely, to ensure the convexity of such $E_{2m}$. By Lempert's work \cite{Le}, all invariant metrics on convex bounded domains, which turn holomorphic maps into Lipschitz $1$-maps, coincide with the Kobayashi metric. Analysis of the non-convex eggs is dealt with in a separate article \cite{BM2}, for clarity.
\medskip\\
\noindent We now provide a description of our results.  First let $m>1$. Then $\partial E_{2m}$ is atleast $C^2$-smooth. Split up $E_{2m}$ as a disjoint union of the following subsets.  
\begin{eqnarray*}
 Z & = & \{ f \big( (0, \hat{0} ) \big) | f \in \mbox{Aut} \; E_{2m} \}, \\
 M^{-} & = & \{ f \big( (p_1, \hat{0} ) \big) | f \in \mbox{Aut} \; E_{2m}, 0 < p_1 < 2^{-1/2m} \}, \\
 M^0 & = & \{ f \big( (2^{-1/2m}, \hat{0} ) \big) | f \in \mbox{Aut} \; E_{2m} \}, \; \mbox{and} \\
 M^+ & = & \{ f \big( (p_1, \hat{0} ) \big) | f \in \mbox{Aut} \; E_{2m}, 2^{-1/2m} < p_1 < 1 \},
\end{eqnarray*}
where $ \mbox{Aut} \; E_{2m} $ denotes the group of holomorphic automorphisms of $ E_{2m} $. Notice that each of the above sets 
is a union of orbits under the action of $ \mbox{Aut} \; E_{2m} $. As orbits of distinct points either coincide or never 
intersect, it follows that the above sets are disjoint. Further, observe that the subdomain $ M^+ $
is a `pinched' one-sided neighbourhood of the strongly pseudoconvex piece of the boundary, pinched along the weakly 
pseudoconvex points of the boundary which form the real algebraic subvariety of 
$\partial E_{2m}$ given by
\[
\omega = \{ (z_1, z_2, \ldots, z_n) : z_1=0, |z_2|^2 + \ldots + |z_n |^2 = 1 \}.
\]
The boundary of $M^+$ is the union of $\partial E_{2m}$ and $M^0$. The significance of these sets is that 
the Wu metric is real analytic, K\"ahlerian of constant negative curvature on $M^+$. All these characteristics break 
down on $M^0$ while just ensuring an overall $C^1$-smoothness of  the Wu metric therein. The real analyticity is 
restored again on $M^{-}$, which lies `inside' $M^0$. Smoothness may be lost along the remaining seam inside $M^0$ namely, $Z$ which however 
encodes some nuances about the smoothness of the boundary near $\omega$, laid down more precisely in the following

\begin{thm} \label{smoothconvexegg}
Let $m $ be any real number bigger than one. The Wu metric on $E_{2m}$ is a $C^1$-smooth Hermitian metric which is real analytic 
on $E_{2m} \setminus (M^0 \cup Z)$. It is of class $C^1$ but not $C^2$ on $M^0$. At points of $Z$, the Wu metric is H\"older of
class $C^{[2m], 2m - [2m]}$ if $m$ is not an integer whereas it is real analytic in case $m$ is an integer. \\
On $M^+$, the Wu metric is K\"ahler of constant negative holomorphic curvature, whereas on $M^-$ it is a non-K\"ahler Hermitian 
metric with its holomorphic curvature non-constant and bounded between a pair 
of negative constants. Overall on $E_{2m}$, the Wu metric has a strongly negative holomorphic curvature in the sense of currents.
\end{thm}

\noindent Note that the set $Z$ is a complex hypersurface in $E_{2m}$ while $M^0$ is a real analytic hypersurface in $E_{2m}$. Indeed, 
$M^0$ is described explicitly by 
\[
\{ z \in E_{2m} \; : \;  2 \vert z_1 \vert^{2m} + \vert z_2 \vert^2 + \ldots + \vert z_n \vert^2 =1 \}
\]
while $Z$ is the intersection of $E_{2m}$ with the hyperplane $\{ z_1 =0 \}$. As the parameter $m$ drops below one, the order of 
tangency of $E_{2m}$ with the normalized ellipsoid namely, the unit ball $ \mathbb{B}^n $, at the boundary orbit accumulation points (which precisely is 
the set of weakly pseudoconvex points in $\partial E_{2m}$) drops below two. The `tangential approach' region $M^+$ to $\omega$, where 
the Wu metric is K\"ahler of constant negative curvature, then disappears altogether.

\begin{thm}\label{convexegg}
Let $1/2 < m <1$. Then $\partial E_{2m}$ is $C^1$ but not $C^2$-smooth. The Wu metric on $E_{2m}$ is a $C^1$-smooth Hermitian metric which is real analytic on $E_{2m} \setminus Z$. At points of $Z$, the Wu metric is $C^{1,2m-1}$ but not of class $C^2$.\\
For $m=1/2$, the Wu metric on $E_{2m}$ is a continuous Hermitian metric, real analytic on $E_{2m} \setminus Z$.\\
For $ 1/2 \leq m < 1 $, the Wu metric on $E_{2m}$ is nowhere K\"ahler. Its holomorphic curvature is bounded above by a negative constant independent of $m$, in the sense of currents. 
\end{thm}

\noindent {\it Acknowledgements:} We thank our advisor Kaushal Verma for suggesting this problem.
\section{Preliminaries}

\noindent Let $ \langle \cdot, \cdot \rangle $ denote the standard Hermitian inner product on $ \mathbb{C}^{n-1} $ which is conjugate-linear in the 
second variable. Write $ z \in \mathbb{C}^n $ as $ z= (z_1, \hat{z})$, where $ \hat{z} = (z_2, \ldots, z_n) $. Any invariant metric, in particular
the Wu and the Kobayashi metric, can be written down explicitly on $E_{2m}$ once it is known at points of $S \cup \{0\}$ where 
 $\{p \in \mathbb{C}^n \; :\; p=(p_1, \hat{0}), \; 0<p_1 < 1 \}$. Indeed, the automorphism 
\begin{equation}\label{aut}
\Phi(z) = \Phi_p(z) = \Big( \frac{|p_1|}{p_1} \frac{s^{1/m}}{ (1-\langle\hat{z},\hat{p} \rangle )^{1/m} } z_1,  \Psi(\hat{z}) \Big)
\end{equation}
of $ E_{2m} $ takes $p=(p_1, \hat{p})$ to $(\tilde{p}_1, \hat{0})$, where $\tilde{p}_1 = |p_1|/s^{1/m}$ with $ s=\sqrt{1- |\hat{p}|^2}$ and
$\Psi$ is the automorphism of $\mathbb{B}^{n-1}$ taking $\hat{p}$ to the origin, given by
\[
\Psi(\hat{z}) = \frac{ \hat{p}I - P(\hat{z}) - sQ(\hat{z}) }{ 1 - \langle\hat{z},\hat{p} \rangle },
\]
where $P$ is the (linear) projection onto the one-dimensional complex subspace spanned by $\hat{p}$, i.e., $P(z) = \frac{\langle\hat{z},\hat{p} \rangle }
{|\hat{p}|^2}{\hat{p}}$ and $Q=I-P$ is the projection  ortho-complementary to $P$. 

\begin{thm} (Theorem 1.6 of \cite{BMV}) \label{T1}
The Kobayashi metric for $E_{2m}$ at the point $(p_1, \hat{0})$ for $0<p_1<1$, is given by
\begin{alignat*}{4} 
K \big( (p_1, 0, \ldots, 0), (v_1, \ldots, v_n) \big) = 
   \left\{ \begin{array}{lrl}
\left( \frac{m^2 p_1^{2m-2} |v_1|^2}{(1-p_1^{2m})^2}
+ \frac{|v_2|^2}{1- p_1^{2m}} + \cdots + \frac{|v_n|^2}{1- p_1^{2m}} \right)^{1/2} & \mbox{for} & u \leq p_1, \\
\\
\frac{m \alpha (1-t) |v_1|}{p_1(1- {\alpha}^2) \left(m(1-t) + t \right) }
& \mbox{for} & u > p_1,
\end{array}
\right.
\end{alignat*}
where
\begin{eqnarray}
u & = & \left( \frac{m^2 |v_1|^2}{|v_2|^2 + \cdots + |v_n|^2} \right)^{1/2}, \label{E4} \\
t & = & \frac{2m^2 p_1^2}{u^2 + 2m(m-1) p_1^2 + u \big( u^2 + 4m(m-1) p_1^2 \big)^{1/2}} \label{E5}
\end{eqnarray} 
and $ \alpha $ is the unique positive solution of 
\begin{equation}
 \alpha^{2m} - t \alpha^{2m-2} - (1-t) p_1^{2m} = 0.
\end{equation}
Moreover, $ K $ is $ C^1$-smooth on $ E_{2m} \times \left( \mathbb{C}^n \setminus \{0\} \right) $ for $ m > 1/2 $. 
\end{thm} 

\noindent The Kobayashi metric at the origin is given by 
$ K ( 0 ,v ) = q_{E_{2m}}(v)$ where $q_{E_{2m}}$ denotes the \textit{Minkowski functional} 
of $E_{2m}$. Overall, the Kobayashi metric is $C^2$ but not $C^3$-smooth, as will be shown in section 6. The first step in investigating the smoothness of the Wu metric is to determine its unit sphere in tangent spaces. In terms of the Euclidean coordinates on the tangent bundle $ E_{2m} \times \mathbb{C}^n $, the unit 
sphere of the Wu metric in $ T_p E_{2m} $ for $p \in S$, is given by
\[
 r_1 |v_1|^2 + r_2 \left( |v_2|^2 + \ldots + |v_n |^2 \right) = 1
\]
where $ r_1 $ and $ r_2 $ are positive real-valued continuous functions of $ p_1 $.  The main objective now is 
\begin{enumerate}
\item[(*)] to find $ r_1, r_2 > 0 $ such that the ellipsoid 
\[
\left\{ v \in \mathbb{C}^n:  r_1 |v_1|^2 + r_2 \left( |v_2|^2 + \ldots + |v_n |^2 \right) < 1 \right\} 
\]
has the smallest volume while containing the set 
\[
I \left( E_{2m}, (p_1, \hat{0}) \right) = 
\left\{ v \in \mathbb{C}^n: K \left( ( p_1, \hat{0} ), v \right) < 1 \right\}.
\]
\end{enumerate}
To this end, first note that the expressions for the best fitting ellipsoid \footnote{The unit ball in any given tangent space with respect 
to the Wu metric will be referred to as the `best fitting ellipsoid' or the Wu ellipsoid.}  at $ (p_1, \hat{0}) $ and the Kobayashi 
indicatrix at $ (p_1, \hat{0}) $ involve terms of the form $ |v_1|^2$, $ |v_2|^2, 
\ldots, |v_n|^2 $ only. Consequently, the problem reduces to finding $ r_1, r_2 > 0 $ such that the set
\[
\left\{ (v_1, \ldots, v_n) : v_1 \geq 0, \ldots, v_n \geq 0, 
r_1 v_1^2 + r_2 \left( v_2^2 + \ldots + v_n^2 \right) = 1 \right\} 
\]
encloses the smallest volume with the coordinate axes and contains the set 
\[
\left\{ (v_1, \ldots, v_n) : v_1 \geq 0, \ldots, v_n \geq 0, 
K \left( ( p_1, \hat{0} ), v \right) \leq 1 \right\}.
\]
\noindent We reformulate problem (*) using the concept of square-convexity from \cite{CK1}. Write 
\begin{eqnarray*}
 x & = &  v_2^2 + \ldots + v_n^2 \; \mbox{and}\\
 y & = &  v_1^2,
\end{eqnarray*}
so that the problem of finding the best fitting ellipsoid at $ (p_1, \hat{0}) \in E_{2m} $ can be restated as:
\begin{enumerate}
 
\item[(**)] Find $ r_1, r_2 $ such that the line segment $ r_1 y + r_2 x = 1 $ in the first quadrant 
bounds the smallest area with the coordinate axes and contains the set 
\[
\left\{ ( x, y ) : v_1 \geq 0, \ldots, v_n \geq 0, 
K \left( ( p_1, \hat{0} ), v \right) \leq 1 \right\}.
\]
\end{enumerate}
The boundary of this set is determined by the equation $ K ^2 \left( ( p_1, \hat{0} ), v \right) = 1 $, which splits into two parts, as is 
evident from Theorem \ref{T1}. Henceforth, the portion of the curve $ K^2 \left( ( p_1, \hat{0} ), v \right) = 1 $ for 
$ u \leq p_1 $ and $ u > p_1 $ will be referred to as the lower K-curve $ C_{low} $ and the upper K-curve $ C_{up} $ respectively. 
In $ x$-$y $ coordinates, the lower K-curve $ C_{low} $ is given by
\begin{equation} \label{E10}
\frac{m^2 p_1^{2m-2} \; y }{(1-p_1^{2m})^2} + \frac{x}{1- p_1^{2m}} = 1,
\end{equation}
which represents a straight line. Denote by $ y = \left( f_{low} ( \sqrt{x}) \right)^2 $, the function representing the lower K-curve $ C_{low} $
in the first quadrant of $ \mathbb{R}^2 $. $C_{low}$ may also be parametrized as:
\begin{eqnarray}\label{E}
x(\alpha) & = & (1-p_1^{2m})^2\alpha, \\
y(\alpha) & =  & (1-p_1^{2m})^2\big( 1- \alpha (1-p_1^{2m})\big)/m^2 p_1^{2m-2} \nonumber
\end{eqnarray}
for $ 1 \leq \alpha < (1-p_1^{2m})^{-1} $. We shall only present the upper K-curve $ C_{up} $ in parametric form as:
\begin{alignat}{3} \label{E8}
x(\alpha) & =  \big( v_2 (\alpha)\big)^2 + \ldots +  \big( v_n (\alpha) \big) ^2 = \frac{ \alpha^{4m-2} + p_1^{4m} - p_1^{2m} \alpha^{2m-2} 
- p_1^{2m} \alpha^{2m} } {\alpha^{4m-2}} \; \mbox{and} \\ \label{E9}
y(\alpha) & =  \big( v_1 (\alpha) \big) ^2 =  \left(\frac{ p_1 \big( m \alpha^{2m-2} - (m-1) \alpha^{2m} - p_1^{2m} \big)}
{m \alpha^{2m-1} } \right)^2, 
\end{alignat}
for $ 0 < \alpha < 1 $ with $ {\alpha}^{2m-1} > p_1^{2m} $. Let $f_{up}$ denote the function representing the upper K-curve, so that $ y( \alpha) = 
\left( f_{up} \big( \sqrt{x(\alpha)} \big) \right)^2 $ is the equation of the upper curve. 

\begin{prop} \label{P1}
The function $ f_{low} $ is both square convex and square concave. The function $ f_{up} $ is strictly square convex for 
$ 1/2 \leq m < 1 $ and strictly square concave for $ m > 1 $. 
\end{prop}

\section{The Wu metric on $ E_{2m} $ for $ 1/2 \leq m < 1 $}
\noindent To find the best fitting ellipsoid for the Wu metric at $p \in S$, first note that
in $ x$-$y $ coordinates, the lower K-curve $ C_{low} $ and the upper K-curve $ C_{up} $ are both convex, by 
Proposition \ref{P1}. Moreover, it follows from $ C^1$-smoothness of the Kobayashi metric on $ E_{2m} $ (refer Theorem \ref{T1}) that 
the upper and the lower K-curve together yield a $C^1$-smooth convex curve in $x$-$y $ coordinates. This implies that the line segment 
obtained by joining the $x$-intercept of the lower K-curve $ C_{low} $ and $y$-intercept of the upper K-curve $ C_{up} $ bounds
the smallest area with the coordinate axes while containing the set $ \left\{ K \left( ( p_1, \hat{0} ), v \right) \leq 1 \right\} $.
It is evident from (\ref{E10}) that the $x$-intercept of the lower K-curve is $ 1 - p_1^{2m} $ while the  
$y$-intercept of the upper K-curve is $ \left( 1- p_1^2 \right)^2 $ since $K  \left( (p_1, \hat{0}), (v_1, \hat{0}) \right) = | v_1|/ \left( 1 - p_1^2 \right) $.
The line segment obtained by joining these two intercepts is given by
\[
\frac{x}{1 - p_1^{2m}} + \frac{y}{ \left( 1 - p_1^2 \right)^2 } = 1.
\]
Observe that the slope of this line segment is greater than the slope of the lower K-curve $ C_{low} $. To summarize, 
$ r_1 = 1/\left( 1- p_1^2 \right)^2  $ and $ r_2 = 1/ \left( 1 - p_1^{2m} \right)$ provides, a solution to the extremal 
problem (**)  and subsequently, the expression for the Wu metric $h_{E_{2m}}$ at reference points:
\begin{multline*}
h_{E_{2m}} (p_1, \hat{0}) = \frac{1} {\left( 1- p_1^2 \right)^2} d z_1 \otimes d \overline{z}_1 + \frac{1} {\left( 1 - p_1^{2m} \right) } 
d z_2 \otimes d \overline{z}_2 + \ldots + \frac{1} {\left( 1 - p_1^{2m} \right) } d z_n \otimes d \overline{z}_n. 
\end{multline*}
\noindent Explicit expression for the Wu metric at $ p \in E_{2m} $, with $ p_1 \neq 0 $ 
is obtained using the automorphisms $\Phi_p$. Continuity of the Wu metric allows its determination at $ ( 0, p_2, \ldots, p_n) \in E_{2m} $.
 \begin{thm} \label{explicit}
For $ 1/2 \leq m < 1 $, the Wu metric on $ E_{2m} $ at the point $ (z_1, \ldots, z_n) $ is given by
\begin{equation*}
 \sum_{i,j=1}^n h_{i\bar{j}} (z_1, \ldots, z_n) \; d z_i \otimes d \overline{z}_j,
\end{equation*}
where 
\begin{alignat*}{3}
 h_{1\bar{1}}  & = \frac{ \left(1 - |\hat{z}|^2 \right)^{1/m} } {  \left( \left(1 - |\hat{z}|^2 \right)^{1/m} - |z_1|^2 \right)^2 }, \\
 h_{1\bar{j}}  & =  \frac{  \left(1 - |\hat{z}|^2 \right)^{-1 + 1/m}  \overline{z}_1 z_j} 
{ m  \left( \left(1 - |\hat{z}|^2 \right)^{1/m} - |z_1|^2 \right)^2 } \; \; \mbox{for} \; \; 2 \leq j \leq n \text{ and } h_{j\bar{1}} = \overline{h_{1\bar{j}}}, \\
h_{j\bar{j}}  & = \left( \frac{  \left(1 - |\hat{z}|^2 \right)^{-2 + 1/m} |z_1|^2 |z_j|^2 } 
{ m^2  \left( \left(1 - |\hat{z}|^2 \right)^{1/m} - |z_1|^2 \right)^2 } + \frac{1 - |\hat{z}|^2 + |z_j|^2 }{ \left( 1 - | \hat{z}|^2 \right) 
\left( 1 - |\hat{z}|^2 - |z_1|^{2m} \right) } \right) \; \mbox{for} \; 2 \leq j \leq n, \\
h_{i\bar{j}}  & = \left( \frac{  \left(1 - |\hat{z}|^2 \right)^{-2 + 1/m} |z_1|^2  \overline{z}_i z_j } 
{ m^2  \left( \left(1 - |\hat{z}|^2 \right)^{1/m} - |z_1|^2 \right)^2 } + \frac{ \overline{z}_i z_j }{ \left( 1 - | \hat{z}|^2 \right) 
\left( 1 - |\hat{z}|^2 - |z_1|^{2m} \right) } \right) \; \mbox{for} \; 2 \leq i, j \leq n \text{ and }  i \neq j .
\end{alignat*}
 \end{thm}
\noindent It follows that the Wu metric is real analytic  on the set $  \{ (z_1, \ldots, z_n ) \in E_{2m} : z_1 \neq 0 \} $. Verifying
$ \partial h_{1 \bar{2}} /\partial z_2 \neq  \partial h_{2 \bar{2}} / \partial z_1$ at points of $S$, by a direct calculation, we infer that the Wu metric is non-K\"ahler. Using Theorem \ref{explicit}, one can also verify that the holomorphic curvature of the Wu metric is bounded between a pair of 
negative constants on the domain $U_Z = E_{2m} \setminus Z$. The strong negativity of the holomorphic curvature current 
across the remaining thin set $Z$ (surrounded by $U_Z$) follows using arguments similar to those outlined in Section 5 (for $ 1/2 < m < 1 $) 
and Section 3 of \cite{BM2} (when $ m =1/2 $). 
%\vspace*{-0.05in}
\begin{rem} \label{degofsmoothness}
Explicit expression of the Wu metric as above, yields information about the overall-order of its smoothness, which relies 
on the smoothness of the power function of one variable $t \to |t|^\alpha$, 
where $\alpha$ is a positive real number. Such a function is infinitely differentiable if $\alpha$ is an integer $>1$, while 
it is only $[\alpha]$-smooth if $\alpha$ is not an integer in which case, the $[\alpha]$th derivative lies in 
the $C^{0, \alpha - [\alpha]}$-H\"older class. We are concerned about the regularity of such power functions 
in a neighbourhood of the origin. Note that $(1-\vert\hat{z} \vert^2)^{1/m}$ is real analytic 
on $E_{2m}$ as it never vanishes there. Barring $m=1/2$, note that the integral part of $2m$ 
is $1$, while its fractional part is $2m-1$. It follows from Theorem \ref{explicit} that the Wu metric 
is \textbf{not} $C^2$-smooth at points of $Z$. However, it is $ C^{1,2m-1}$-smooth at $ Z $. So the 
Wu metric is $C^1$-smooth overall, unless $m=1/2$.
\end{rem}
\section{The Wu metric on $ M^+ $ for $ m > 1 $}
\begin{thm}
The Wu metric tensor at $p=(p_1,p_2, \ldots,p_n) \in M^{+}$ is given by
\begin{multline*}
\frac{1}{(s^2-|p_1|^{2m})^2} \Big( m^2s^2|p_1|^{2m-2} dp_1 \otimes  d \overline{p}_1 + \sum\limits_{j=2}^{n}(s^2+|p_j|^2-|p_1|^{2m}) dp_j \otimes  d \overline{p}_j \\
+ 2 \Re \big( \sum\limits_{j=2}^{n} m|p_1|^{2m-2} p_1  \overline{p}_j dp_j \otimes  d \overline{p}_1 \big)
+ 2 \Re \big( \sum\limits_{\substack{{j,k=2}\\ {k<j}}}^{n} p_k \overline{p}_j dp_j \otimes d \overline{p}_k \big) \Big) 
\end{multline*}
\end{thm}
\noindent The smooth function $\rho(z) = - \log\big( 1- (\vert z_1 \vert^{2m} + \vert z_2 \vert^2 + \ldots + \vert z_n \vert^2 ) \big)$ is a K\"ahler potential for the Wu metric on $M^{+}$. It is also straightforward to verify that the Wu metric is  K\"ahler  with constant holomorphic curvature $-2$ on $M^+$.
\section{The Wu metric on $ M^- $ for $ m > 1 $}
\begin{prop} \label{P4}
For $ 0 < p_1 < 2^{-1/2m} $, the unit sphere of the Wu metric in $ T_{(p_1, \hat{0})} E_{2m} $ does not intersect the 
lower K-curve $ C_{low} $ and intersects the 
upper K-curve $ C_{up} $ in a unique point $ (x^*, y^*) $ with $ x^* \neq 0 $.
\end{prop}
\noindent Real analytic dependence of $ (x^*, y^*) $ on the base point, can be deduced by the arguments involving  Fritz John's generalization of the Lagrange multiplier method to study the variation of the best fitting 
ellipsoid (as we move from one tangent space to another), as in \cite{CK1}. The `Lagrangian' functional adapted to 
study all the constraints and parameters jointly, is of the form  
\begin{equation} \label{E0}
H=H(p, R_1(p), R_2(p), V_1(p), V_2(p), \Theta(p) \big) : \big( 0, 2^{-1/2m} \big) \times ( 0, \infty)^5 \rightarrow \mathbb{R}^5
\end{equation}
for our case as well. This will be useful for analyzing 
the regularity of the Wu metric across the thin sets $M^0$ and $Z$ at the boundary of $M^{-}$. For the analysis of curvature, it is enough 
to determine the curvature tensor at points of the form $(p_1 ,\hat{0})$ where the matrix representing the Wu metric is diagonal. This 
requires knowledge about the metric tensor in a neighbourhood of points $(p_1 ,\hat{0})$ which we record below. First, let $\alpha^*$ 
be the value of the parameter $ \alpha $, where the upper $K$-curve meets the square transform 
of the best fitting ellipsoid. Writing $X= (\alpha^*)^2$, we remark that $X$ may also be defined as the unique solution in the interval $(0,1)$ of the equation
\[
s^4 X^{2m-1} -(m+1)|p_1|^{2m}s^2X^{m-1} + (m-2)|p_1|^{2m}s^2X^m + 2|p_1|^{4m} =0.
\]

\begin{prop} \label{M-exp}
\noindent  The Wu metric at $p \in M^{-}$ is given by
\begin{multline*}
\frac{s^2X^{2m-1}}{2F_s^2} \Big( \frac{m^2s^2}{|p_1|^2}  dp_1 \otimes  d \overline{p}_1 + \sum\limits_{j=2}^{n} \big(\frac{F_{s_j}}{|p_1|^{2m}} -1 \big) dp_j \otimes  d \overline{p}_j
+2 \Re \sum\limits_{j=2}^{n}\frac{m\overline{p}_j}{\overline{p}_1} dp_j \otimes  d \overline{p}_1\\
+\frac{mX^{m-1}-(m-1)X^m}{|p_1|^{2m}}\sum\limits_{j,k=2, j\neq k}^{n} \overline{p}_j p_k dp_j \otimes  d \overline{p}_k \Big)
\end{multline*}
where $F_s = ms^2X^{m-1} - (m-1)s^2X^m - |p_1|^{2m}$, $s_j^2= s^2 + |p_j|^2 = 1- (|\hat{p}|^2 - |p_j|^2)$ 
and $F_{s_j}= ms_j^2X^{m-1} - (m-1)s_j^2X^m$. Moreover, the Wu metric is real analytic on the set $ M^- $.
\end{prop}
\noindent It turns out that on $M^{-}$, each component of the curvature tensor either vanishes identically or is a strictly negative 
real valued function. We now list the components $R_{i \overline{j} k \overline{l}}$ of the curvature tensor, with 
(standard) notations as in \cite{CK1}, which are negative. First, consider 
the case where atmost two of the indices $(i,j,k,l)$ are distinct and one of them is $1$. Let $\gamma \in  \{2, 3, \ldots,n\}$. It follows from \cite{CK1} 
that if $(i,j,k,l)$ is one of 
\[
(\gamma, \gamma, \gamma,\gamma), (\gamma, \gamma, 1,1),(\gamma,1,1, \gamma), (1,\gamma, \gamma, 1), (1,1,\gamma, \gamma), (1,1,1,1)
\]
then $R_{i \overline{j} k \overline{l}}$ is strictly negative. For quadruples of $1, \gamma$ not listed above, the corresponding curvature component vanishes. It can 
be checked that $R_{1 \overline{1} \gamma \overline{\gamma}} \neq R_{\gamma \overline{1} 1 \overline{\gamma}}$ at $(p_1, \hat{0})$ for all $p_1 \in (0, 2^{-1/2m})$, and
hence the Wu metric is nowhere K\"ahler on the orbit of the segment $(0, 2^{-1/2m})$, i.e. on $M^{-}$. Next, if $i=j>1$ and $k=l>1$, then  
$ R_{i \overline{i} k \overline{k}}$ is negative. In case, both the indices $i,j$ are distinct and  $>1$, it turns out that  $R_{i \overline{j} k \overline{l}}$ is negative, iff  $(k,l)= (j,i)$. In \textbf{all} other cases, $R_{i \overline{j} k \overline{l}}=0$. It is then possible to establish
assertions about the curvature on $M^{-}$ as in Theorem 1.2.
\section{The Wu metric along and across the thin sets $Z$ and $M^0$}
\noindent Analysis of the Wu metric in the remaining thin subsets $Z$ and $M^0$ of $E_{2m}$, relies firstly on 
the smoothness of elementary power functions, as mentioned in remark \ref{degofsmoothness}. We shall deal only with the case $m \not \in \mathbb{N}$, as 
the complementary case is easier. Let us begin by analyzing the smoothness of the Wu metric {\it at} the origin. This entails studying the behaviour of how the Kobayashi metric is changing at the origin. We recall that 
\[
K \big((0, 0, \ldots, 0), (v_1, \ldots,v_n) \big)= \frac{1}{\tilde{\alpha}},
\]
where $\tilde{\alpha}$ is the unique positive solution of the equation
\[
\vert v_1 \vert^{2m} \tilde{\alpha}^{2m} + \big( \vert v_2 \vert^2 + \ldots + \vert v_n \vert^2 \big) \tilde{\alpha}^2 =1. 
\]
It follows that $\tilde{\alpha}$ is $C^{[2m]}$-smooth as a function of $\vert v_1 \vert$, and a real analytic function of the variables $ \vert v_2 
\vert, \ldots, \vert v_n \vert$. Hence, $\tilde{\alpha}$ is $C^{[2m]}$-smooth function of the variables $ \vert v_1 \vert, \vert v_2 \vert, \ldots, \vert v_n \vert$. Indeed, these conclusions can be derived by consulting the assertions concerning the regularity
of the solution of $C^k$-smooth equations, in the implicit function theorem. We now examine the regularity of $K(p,v)$, when $p$ varies across $Z$. It suffices 
to focus attention on the case when $p$ varies through $S^\epsilon = \{(p_1, \hat{0}):  p_1 \in [0,\epsilon)\}$ for some $\epsilon>0$. We shall 
simultaneously deal with joint regularity in the variables $p,v$ as well, though we shall be terse; in this regard, observe that the 
Kobayashi metric at a point $p$ in a small neighbourhood of the origin, is of the form $K(\Phi_p(p), D\Phi_p(v))$ with $\Phi_p(z)$ as in 
(\ref{aut}). Note that $\Phi_p(z)$ is jointly real analytic in the variables $z$ and (the parameter) $p$. Consequently, the $C^{[2m]}$-smoothness of $K(\cdot,\cdot)$
in a neighbourhood of the origin follows, as soon as we verify that $K(p,v)$ is $C^{[2m]}$-smooth for $p$ varying in $S^\epsilon$. As 
for $v$, it suffices to restrict attention in a neighbourhood of the set of points of contact of the Wu ellipsoid with the Kobayashi indicatrix
in $T_pE_{2m}$. Now, recall from Proposition \ref{P4} that, on $M^{-}$, the contact point of the square transforms of the Kobayashi indicatrix 
and the Wu ellipsoid, lies on the upper K-curve, where the 
Kobayashi metric is described as follows (cf. equation (7.29) of \cite{BMV}):
\begin{equation} \label{Kobalternate}
\Big(K \big( (p_1,\hat{0}), v \big) \Big)^2 = \frac{\tilde{\alpha}^2 \vert v_1 \vert^4}{\tilde{t}^2\Big( 
\vert v_1\vert^2 \tilde{\alpha}^2 - p_1^2 \tilde{\alpha}^{2m} \vert v_1 \vert^{2m} \Big)^2 }.
\end{equation}
Here, $x = \tilde{\alpha}$ satisfies the equation $F(v,p,x)=0$ with $p=(p_1, \hat{0})$, where
\begin{equation} \label{eqnsol}
F(v,p,x) := \Big( 1- \tilde{t}^2 \frac{\vert \hat{v} \vert^2 p_1^2}{\vert v_1 \vert^2}\Big) \vert v_1 \vert^{2m} x^{2m} 
+ \tilde{t}^2 \vert \hat{v} \vert^2 x^2 - 1,
\end{equation}
and
\begin{equation} \label{teqn}
\tilde{t}^2 = \frac{2 \vert v_1 \vert^2}{\vert v_1 \vert^2 + 2 \left(1-\frac{1}{m} \right) \vert \hat{v} \vert^2 p_1^2 + \vert v_1 \vert 
\sqrt{\vert v_1 \vert^2 + 4 \left(1-\frac{1}{m} \right) p_1^2 \vert \hat{v} \vert^2}}.
\end{equation}
An application of the implicit function theorem to $F(v,p,x)$  will enable us to
write $\tilde{\alpha}$ as $ C^{[2m]} $-smooth function of $p$ and $v$, provided we verify 
$
\frac{\partial F}{\partial x}(w,p, \tilde{\alpha}) \neq 0,
$
where $ w = (w_1, w_2, \ldots, w_n) $ is any of the points corresponding to the point $(x^*, y^*)$, obtained in Proposition \ref{P4}; so,
$ w_1^2 = y^*$ and $ w_2^2 + \ldots + w_n^2 = x^* $. To obtain a contradiction, assume that $\partial F/\partial x (w,p, \tilde{\alpha})=0$ i.e.,
\[
\frac{\partial F}{\partial x} (w,p, \tilde{\alpha}) = 2m \Big( 1-  \tilde{t}^2 \frac{\vert \hat{w} \vert^2 p_1^2}
{\vert w_1 \vert^2  } \Big) \vert w_1 \vert^{2m} x^{2m-1} + 2 \vert \hat{w} \vert^2 \tilde{t}^2 x = 0.
\]
 Since $ \tilde{\alpha} \neq 0$, it follows from the 
above equation that 
\begin{equation}\label{alphre}
\Big( 1-  \tilde{t}^2 \frac{ \vert \hat{w} \vert^2  p_1^2}{\vert w_1 \vert^2} \Big)\vert w_1 \vert^{2m} {\tilde{\alpha}}^{2m} = - \frac{ \vert \hat{w} 
\vert^2 \tilde{t}^2 {\tilde{\alpha}}^2} {m},
\end{equation}
so that the defining equation for $\tilde{\alpha}$ can be rewritten as 
$
\tilde{\alpha}^{2} = m(m-1)^{-1}\vert \tilde{t} \hat{w}\vert^{-2}.
$
Substituting the above expression for $ \tilde{\alpha} $ into (\ref{alphre}), dividing throughout by $ \tilde{t}^2 \vert \hat{w} \vert^2 p_1^2$ yields
\[
\Big( \frac{1}{p_1^2 \tilde{t}^2} \frac{\vert w_1 \vert^2}{\vert \hat{w} \vert^2 } - 1 \Big) \vert w_1 \vert^{2m-2} = - 
\frac{(m-1)^{m-1}}{m^m}\frac{ \big( p_1^2\tilde{t}^2\vert \hat{w} \vert^2  \big)^{m-1} }{ p_1^2} \big( \frac{1}{p_1^2} \big)^{m-1}
\]
Define $
f(p_1,w) = \vert w_1 \vert^2/p_1^2\tilde{t}^2\vert \hat{w} \vert^2$,
so as to write the last equation as
\begin{equation} \label{freritten}
\big( f(p_1,w) - 1 \big) \big( f(p_1,w) \big)^{m-1} = - \frac{(m-1)^{m-1}}{m^m} \frac{1}{p_1^{2m}}.
\end{equation}
We claim that $f(p_1,w)>1$. To prove this claim, note that
\[
f(p_1,w) = \frac{R^2}{2 p_1^2} + \frac{R}{2 p_1} \Big( \frac{R^2}{p_1^2} + 4 \big(1- \frac{1}{m} \big) \Big)^{1/2} 
+ 1- \frac{1}{m},
\]
where $R$ denotes the ratio $\vert w_1 \vert / \vert \hat{w} \vert$.
At this point, recall that we are in the case $u \geq p_1$ (where $u = m \vert w_1 \vert / \vert \hat{w} \vert$) or equivalently, $R/p_1 > 1/m$. Using 
this in the above equation gives the claim. This is a contradiction since the right hand side of (\ref{freritten}) is 
negative. Hence, we conclude that our assumption $\frac{\partial F}{\partial x}(w,p, \tilde{\alpha}) =0$  must be
wrong. This enables an application of the implicit function theorem to deduce that $\tilde{\alpha}$ is a $[2m]$-smooth function of $ (p, v) $. Next, (\ref{teqn})  
shows that $\tilde{t}$ is a smooth function. Further, notice that the expression (\ref{Kobalternate}) of the Kobayashi metric is a rational function 
(whose denominator is well-defined on $E_{2m}$) involving elementary power functions of $v, p, \alpha,\tilde{t}$. Hence, $ K^2 $ is $C^{[2m]}$-smooth 
function of $(p, v) $ near the origin. Next, another application of the implicit function theorem, this time to the functional $H$, 
mentioned at (\ref{E0}), one gets the desired smoothness of the Wu metric near $Z$. It may seem that, in doing so, there may be a loss of a degree of smoothness for the Wu metric as compared to the $C^{[2m]}$-smoothness of the Kobayashi metric. However, it is
possible to discern the optimal order of smoothness of the Wu metric to be $C^{[2m]}$ (as in Theorem 4.1 of \cite{Kim}), owing
to the fact that that the meeting point of the Wu ellipse with the $K$-curve in the tangent spaces near the origin, lies away from the
coordinate axes; this ensures the real analytic dependence of $K^2$ on the variables $v$. This together with an analysis of the explicit 
expression of $K^2$ yields the $C^{[2m]}$-smoothness of the component functions of $H$; thereby the $C^{[2m]}$-smoothness of the Wu metric near the origin, subsequently, near $Z$. To remark now about the case $m \in \mathbb{N}$, similar analysis yields the conclusion that the Wu metric is real analytic near $Z$ in this case.  
\medskip\\
\noindent In order to analyse the Wu metric on $M^0$, we do {\it not} need to assume that $m \not \in \mathbb{N}$ for convenience. We only
provide a couple of details to support the arguments in \cite{CK1}, for passing through to our present setting. Firstly, we claim that the Kobayashi
indicatrix is $C^2$-smooth but not $C^3$-smooth on $ E_{2m} $, for $m>1$. To establish this, it suffices to analyze the smoothness at the 
joining point $ u = p_1 $ of the upper and lower K-curves. Note that the lower K-curve (as described by equation (\ref{E10})) is a straight line and hence 
$ \partial^l y/\partial x^l \equiv 0$ for $l \geq 2$. We now show that, for the upper K-curve $ \big( x(\alpha), y(\alpha) \big)$, $\partial^l y/\partial x^l$ 
is zero for $l=1,2$ but not for $l=3$ at the joining point. To this end, note that the 
joining point $u=p_1$ corresponds to the limiting value $1$ for the parameter $\alpha$. 
\begin{equation}\label{2ndder}
\frac{\partial^2 y}{\partial x^2} =\frac{\dot{x}(\alpha)\ddot{y}(\alpha) - \dot{y}(\alpha) \ddot{x}(\alpha)}{\big(\dot{x}(\alpha) \big)^3}.
\end{equation}
A direct computation shows that the numerator of the right hand side above is given by
\[
-  8m(m-1)\alpha^{2m-3}(1- \alpha^2)p_1^{2m+2} \big( m \alpha^{2m-2} +(m-1)\alpha^{2m}-(2m-1)p_1^{2m} \big)^2/m^2\alpha^{8m-2},
\]
which as a function of $p_1$ vanishes identically when $ \alpha = 1 $.
This confirms that the values of $\partial^2 y/\partial x^2$ of the upper and lower K-curves
match at the ``joining point'' and leads to the $C^2$-smoothness of the Kobayashi indicatrix in $T_pE_{2m}$ with $p=(p_1, \hat{0})$. 
Next, we evaluate $\partial^3 y/\partial x^3$  and show that the values fail to match up. Indeed, as the 
numerator in (\ref{2ndder}) is zero at $\alpha =1$, it can be calculated that at the joining point, $\partial^3 y/\partial x^3$ is given by
\begin{equation}\label{3rdder}
\big( \dot{x}(1) \dddot{y}(1) - \dot{y}(1) \dddot{x}(1)  \big)/(\dot{x}(1))^4
\end{equation}
provided $\dot{x}(1) \neq 0$ -- it turns out that $\dot{x}(1) \neq 0$ precisely when $m=1/2$. The numerator of (\ref{3rdder}) is 
$ 16p_1^{2m+2} (1-p_1^{2m})^2(2m-1)^2(m-1)/m$. As $1-p_1^{2m}$ never vanishes on $E_{2m}$ and we are in the case $m \neq 1$ and more importantly $m \neq 1/2$, this
proves that the Kobayashi metric is not $C^3$ smooth on $E_{2m}$ when $m>1$. For the case $m<1$, we only need to look at the indicatrix at 
the origin, to conclude that the Kobayashi metric is not of class $C^2$. Getting back again to the case $m>1$, it remains to study how these indicatrices vary with respect to $ p $,
to establish $C^2$-smoothness of the Kobayashi metric. Indeed, first note the real analyticity in the expression for the `joining point' given by: $ \left((1-p_1^{2m})^2, p_1^2(1-p_1^{2m})^2/m^2 \right)$.
While this provides a clue about the $C^2$-smoothness, for a proof it remains to check that $\partial^2/\partial p_1^2$ and 
the $\partial^2/\partial p_1 \partial x$ derivatives of $(f_{up}(\sqrt{x}))^2$ and $(f_{low}(\sqrt{x}))^2$ (which are actually functions of both $x$ and $p_1$)
agree when evaluated at the joining point. This indeed happens and involves several implicit differentiations. For this we note that (\ref{E}), (\ref{E8})  do {\it not} together form a parametrization of the $K$-curve, which is of desired smoothness at the joining point. Nevertheless, a  calculus with these parametrizations confirms that the $K$-curves do coalesce to form a regular surface in the $(x,y,p)$-space, which is $C^2$-smooth. More precisely, {\it all} second order partial derivatives 
of $(f_{up}(\sqrt{x}),p_1)^2$ and $(f_{low}(\sqrt{x},p_1))^2$ match as required. This  leads to the $C^2$-smoothness of the the aggregate of the indicatrices, thereby yielding  $C^2$-smoothness of the Kobayashi metric, first at reference points of the form $(p_1, \hat{0})$ for $0<p_1<1$; thereafter, also on $E_{2m}\setminus Z$ which is the orbit of the segment $S$ under the real analytic action 
of $\textrm{Aut}(E_{2m})$. Since we have already discussed the smoothness of the Kobayashi metric in
a neighbourhood of the origin, together with the foregoing observations,  this completes the verification that the optimal overall smoothness of the Kobayashi metric is indeed $C^2$, when $ m > 1 $. One may now proceed with  arguments employed in \cite{CK1} (past Proposition 7 therein) to conclude the regularity-analysis on $M^0$.\medskip\\
\noindent We now conclude the article with remarks, for completeness, on the curvature of the Wu metric on the critical sets $M^0$ and $Z$. It is enough to focus attention near the point $p=(2^{-1/2m},0)$ in $M^0$. As $M^0$ is a finite type hypersurface, it cannot contain any non-trivial germ of a complex analytic variety. In particular, any Riemann surface through $p$ can intersect $M^0$ only at $p$ or along a smooth curve. This together with the fact that 
the holomorphic curvature is bounded above by a negative constant on either sides of $M^0$ namely, $M^-$ and $M^+$, paves the way for 
applying the  lemma in Appendix B of \cite{CK1}, thereby rendering the negativity of the holomorphic curvature current on $M^0$. Recall that $ Z $ is 
a complex hypersurface and consequently, $ Z $ intersects a embedded Riemann surface in $E_{2m}$ only in 
a discrete set of points, unless it is contained in (as an open subset of) $Z$. In the former case, the required negativity of the
holomorphic curvature current follows from the arguments using the lemma in Appendix B of \cite{CK1}, as the holomorphic curvature
is strongly negative in any small (deleted) neighbourhood of $Z$ consisting of points of $M^-$. In the latter case, we may assume 
that the Riemann surface is an open subset $U$ of $Z$ containing the origin. Pick any one-dimensional complex subspace $l$ in $Z$ -- recall 
that $Z$ is in fact a complex hyperplane. So $ l_U =l \cap U $ is contained in $l$. We restrict attention to a small disc $ D $ about the origin
with $D \subset l_U$. Note that the holomorphic curvature in the direction of $l$ is realized by the Gaussian curvature of the Wu metric
restricted to $D$. As the Wu metric is an invariant metric, the desired negativity of the holomorphic curvature follows.

\end{document}